\documentclass[a4paper,12pt]{amsart}
\usepackage[cp1251]{inputenc}
\usepackage[russian]{babel}
\usepackage{amssymb, amsmath, latexsym}
\usepackage{amsfonts,mathtext,cite,enumerate,float,amsthm}
\numberwithin{equation}{section}
\renewcommand{\le}{\leqslant}
\renewcommand{\ge}{\geqslant}

\newenvironment{Proof}
{\par\indent{\bf Доказательство.}}
{\hfill$\scriptstyle\blacksquare$}
\newtheorem{theorem}{\indent Теорема}
\newtheorem*{thS-L}{\indent Теорема Сильвестра-Любича}

\newtheorem{lemma}{\indent Лемма}

\theoremstyle{definition}
\newtheorem*{definition}{\indent Определение}
\newtheorem{example}{\indent Пример}
\newtheorem{remark}{\indent Примечание}

\usepackage{geometry} 
\begin{document}

\begin{center}

\bf

{\large

Extraction of harmonics\\ from trigonometric polynomials\\ by amplitude and phase operators

}

\bigskip

V.~I.~Danchenko and D.~G.~Vasilchenkova

\end{center}

\bigskip

Extraction of harmonics of a given order from real trigonometric polynomials
(signals) is one of the main problems in harmonic analysis.
It has many applications in physics, radio and electrical engineering,
in particular, in filtration of harmonic signals of different nature.
There exist many methods (mainly, approximative) for solution of this problem.
The most common ones are spectral methods based on Fourier transform and other resonance principles.

In this paper we propose a new method for extracting harmonics by amplitude and phase transformation
of trigonometric polynomials. These transformations use the two simplest operations --- multiplication by a real constant and phase shift  ---
to obtain polynomials \textit{similar} to the initial ones. A harmonic is extracted by an amplitude and phase operator
that simply overlays (sums up) a finite number of such \textit{similar} polynomials. The overlay method enables us to obtain precise
analytical formulas for calculating harmonics of a given order.

\bigskip


\newpage

\title[Выделение гармоник аплитудно-фазовыми операторами]
{Выделение гармоник \\ из тригонометрических многочленов \\
аплитудно-фазовыми операторами}
\author{В.\,И.~Данченко, Д.\,Г.~Васильченкова}

\footnotetext{Работа выполнена при финансовой поддержке ДРПННиТ
$N^{0}$ 1.1348.2011, РФФИ (проект $N^{0}$ 14-01-00510),
Минобрнауки России (задание $N^{0}$ 1.574.2016/ФПМ).}


\maketitle


\begin{abstract}
Выделение гармоник заданного порядка из вещественных
тригонометрических полиномов (сигналов) --- одна из важных задач
гармонического анализа. Она имеет множество приложений в физике,
радиотехнике, электротехнике, например, в процессах фильтрации
гармонических сигналов разной природы. Разработано немало методов
(в основном приближенных) решения этой задачи. Наиболее известными
являются спектральные методы на основе Фурье-преобразований и
других резонансных принципов.

В настоящей работе для выделения гармоник предлагается новый метод
аплитудно-фазовых преобразований тригонометрических многочленов.
Они переводят многочлены в {\it подобные} им, совершая две
простейшие операции --- домножение на вещественную константу и
сдвиг по фазе. Гармоника выделяется амплитудно-фазовым оператором,
представляющим собой простое наложение (сложение) конечного числа
подобных многочленов. Метод наложения позволил получить точные
аналитические формулы для вычисления гармоник заданного порядка.
\end{abstract}

\section{Введение. Постановка задачи}\label{s1}

\subsection{}\label{s1.1}

Пусть имеется тригонометрический многочлен
$$
T_n(x)=\sum_{k=1}^{n}\tau_k(x),\quad \tau_k(x):=a_k\cos kx+b_k\sin
kx,\qquad a_k,\;b_k\in {\mathbb R},\quad n\in {\mathbb N}
$$
(для удобства мы рассматриваем многочлены без свободного члена).
Амп\-ли\-туд\-но-фазовым оператором (АФО) порядка $\le n$ этого
многочлена будем называть преобразование вида
\begin{equation}
H_m(T_n,\{X_j\},\{\lambda_j\};x)=\sum_{j=1}^{n} X_j\cdot
T_n\left(x-\lambda_j\right), \qquad X_j,\,\lambda_j\in{\mathbb R}
  \label{(1)}
\end{equation}
(порядок АФО равен числу различных слагаемых (\ref{(1)}), в
которых $\exp(i\lambda_j)$ попарно различны, а $X_j\ne 0$). Наша
основная

{\bf Задача.} {\it Подбором амплитуд $X_j$ и начальных фаз
$\lambda_j$ получить гармонику с заданным номером
${\mu}=\overline{1,n}$ в виде}:
\begin{equation}
\tau_{\mu}(x)=H_n(T_n,\{X_j\},\{\lambda_j\};x).
 \label{(2)}
\end{equation}

\begin{definition}
Решение $X_j$, $\lambda_j$ будем называть допустимым, если
выполнено условие $(\mathcal A)$:\;\textit{параметры $X_j$,
$\lambda_j$ вещественны и не зависят от $T_n$.}
\end{definition}
Отметим, что если некоторые $X_j= 0$, или не все
$\exp(i\lambda_j)$ попарно различны, то порядок АФО $(\ref{(1)})$
будет, очевидно, строго меньше $n$.
\begin{definition} Решение $X_j$, $\lambda_j$ будем называть
регулярным, если выполнено условие $({\mathcal A}')$:\; \textit{
порядок АФО в точности равен $n$, т.е. если в $(\ref{(1)})$ все
$\exp(i\lambda_j)$ попарно различны, а $X_j\ne 0$.}
\end{definition}

\begin{remark} \label{remark0}

Отметим, что основная трудность задачи (\ref{(2)}) состоит в том,
что АФО должен быть {\it вещественным и иметь порядок} $\le n$.
Действие вещественного АФО имеет простой физический смысл
фильтрации: из смеси $T_n$ гармоник выделяется гармоника
$\tau_{\mu}(x)$ наложением  подобных сигналов, отличающихся лишь
амплитудами и начальными фазами.

Задача выделения отдельных гармоник из сигнала --- одна из важных
задач гармонического и спектрального анализа. Она имеет
многочисленные приложения в физике, радиотехнике, электротехнике,
например, в процессах фильтрации гармонических сигналов разной
природы (см., \cite{Yazdekhasti} ---\cite{SIBERT}).

\end{remark}

\subsection{}\label{s1.2}
Кратко остановимся на основных методах и результатах работы.
Несложно показать (см. п. \ref{s1.2}), что задача (\ref{(2)})
равносильна неполной системе дискретных моментов с неизвестными
$Y_j:=X_j z_j$ и $z_j:=e^{-i\lambda_j}$:
\begin{equation}
\left\{
\begin {array}{ccc}
Y_1+Y_2+\ldots+ Y_{n} &= &{\sigma}_1\\
\ldots & \ldots &\ldots\\
Y_1 z_1^{\mu-1}+Y_2 z_2^{\mu-1}+
\ldots+ Y_{n} z_{n}^{\mu-1} &= &{\sigma}_{\mu}\\
\ldots & \ldots &\ldots\\
 Y_1 z_1^{n-1}+Y_2 z_2^{n-1}+\ldots+ Y_{n} z_{n}^{n-1} &= &{\sigma}_n
\end {array}
\right.,\quad Y_j=X_j z_j,\quad z_j=e^{-i\lambda_j},
 \label{(4)}
\end{equation}
при условии $(\mathcal A)$, которое можно записать в другой форме:

\smallskip \textbf {$(\mathcal A)$}: \textit {${\sigma}_k=0$ при $k\ne {\mu}$ и
${\sigma}_{\mu}=1$, $z_j=e^{-i\lambda_j}$, все $\lambda_k$ и $X_k$
вещественны $($некоторые $X_k$ могут быть нулевыми$)$.}

Как будет видно из дальнейшего (см. раздел \ref{s2}), при решении
задач (\ref{(2)}) и (\ref{(4)}) основную роль играет
вспомогательный вещественный параметр
\begin{equation}
\omega:=X_1+X_2+\ldots+ X_{n},
 \label{(001)}
\end{equation}
подбором которого и удовлетворяется условие $(\mathcal A)$.
Наиболее полные результаты получены при дополнительном условии:
\begin{equation}
n=s\mu-1,\qquad s\in {\mathbb N}.
 \label{(000)}
\end{equation}
В этом случае получены явные формулы для решений задач (\ref{(2)})
и (\ref{(4)}). Основная идея решения состоит в следующем. Система
(\ref{(4)}) дополняется до полной системы моментов (см.
(\ref{(15)})) добавлением уравнения (\ref{(001)}) и уравнений,
комплексно сопряженных с (\ref{(4)}). Полученная система
(\ref{(15)}) в силу условий $(\mathcal A)$ равносильна системе
(\ref{(4)}). Важно то, что вещественность амплитуд $X_k$ в
совместной системе (\ref{(15)}) выполняется автоматически (лемма
\ref{lemma1}). Для решения полных систем моментов обычно
применяется классический метод Прони \cite{Prony}, причем для его
эффективности требуется дополнительно условие регулярности
$({\mathcal A}')$. Однако при $\mu\ge 2$ и условии (\ref{(000)})
все решения задач (\ref{(2)}) и (\ref{(4)}) не являются
регулярными (теорема \ref{theorem2}), что приводит к невозможности
применения метода Прони. Для таких случаев в разделе \ref{s4}
разработано несколько методов регуляризации задачи (\ref{(15)})
определенными вариациями правых частей. Это позволило предельными
переходами из решений регуляризованных систем получить допустимые
решения задач (\ref{(2)}) и (\ref{(4)}) (теоремы \ref{theorem6},
\ref{theorem7}).

Методами регуляризации показано, что для разрешимости задач
(\ref{(2)}) и (\ref{(4)}) достаточно, чтобы параметр $\omega$
являлся корнем уравнения $U_{s}\left({\omega}/{2}\right)=0$, где
$U_{s}$ --- многочлен Чебышева второго рода степени $s$.

В случаях $\mu=1,2$ это условие является и необходимым (теоремы
\ref{theorem3},\ref{theorem6}), что позволяет найти {\it все}
допустимые и регулярные решения $\{X_j(\mu,\omega)\}$,
$\{\lambda_j(\mu,\omega)\}$.

В каждом случае $\mu\ge 1$ при условии (\ref{(000)}) имеет место
неединственнось: существует $s$ различных решений
$\{X_j(\mu,\omega_k)\}$, $\{\lambda_j(\mu,\omega_k)\}$,
соответствующих $s$ различным корням $\omega_k$ уравнения
$U_{s}\left({\omega}/{2}\right)=0$ (теоремы \ref{theorem3},
\ref{theorem6}, \ref{theorem7}).

Отметим, что для удобства мы рассматриваем многочлены без
свободного члена $a_0$. В общем случае, когда $a_0\ne 0$, вместо
(\ref{(2)}) имеем
\begin{equation}
a_0\,\omega+\tau_{\mu}(x)=H_n(T_n,\{X_j(\omega)\},\{\lambda_j(\omega)\};x),
 \label{(2')}
\end{equation}
где $\omega$ --- корень уравнения
$U_{s}\left({\omega}/{2}\right)=0$.

Отметим еще, что для найденных решений равенство (\ref{(2')})
выполняется не только на многочленах, но и на достаточно широком
классе сходящихся тригонометрических рядов (см. примечания
\ref{remark1a}, \ref{remark2}). Например, в случае выделения
первой гармоники равенство (\ref{(2')}) выполняется на
тригонометрических рядах, в которых отсутствуют гармоники с
номерами $\pm 1+(n+2)k$, $(n+2)k$, $k\in{\mathbb N}$.

\section{Система уравнений для амплитуд и фаз}\label{sss1}

\subsection{}\label{s1.2}
Имеем
$$
\tau_k(x)
 =\frac{a_k}{2}\left(e^{ixk}+e^{-ixk}\right)-
 i\frac{b_k}{2}\left(e^{ixk}-e^{-ixk}\right)=
$$
$$
 =\left(\frac{a_k}{2}-i\frac{b_k}{2}\right)e^{ixk}+
\left(\frac{a_k}{2}+i\frac{b_k}{2}\right)e^{-ixk}=
 \alpha_k e^{ixk}+\overline{\alpha_k e^{ixk}}.
$$
Поэтому при $z_j=e^{-i\lambda_j}$ получаем
$$
T_n\left(x-\lambda_j\right)=\sum_{k=1}^{n} \tau_k(x-\lambda_j)=
$$
$$
 =\sum_{k=1}^{n} \left(\alpha_k e^{ik\left(x-\lambda_j\right)}+
 \overline{\alpha_k e^{ik\left(x-\lambda_j\right)}}\right)
 =2{\rm Re}\sum_{k=1}^{n} \alpha_k z_j^k e^{ik x}.
$$
Следовательно, при вещественных $X_j$ имеем
\begin{equation}
H_n(x)=\sum_{j=1}^{n} X_j\cdot T_n\left(x-\lambda_j\right)
 =
2\cdot{\rm Re}\left(\sum_{k=1}^{n} \left(\sum_{j=1}^{n} X_j
z_j^k\right) \alpha_k e^{ik x}\right).
  \label{(3)}
\end{equation}
Таким образом, задача (\ref{(2)}) равносильна системе моментов
(\ref{(4)}) с условием $(\mathcal A)$. Действительно, если решение
(\ref{(4)}) найдено, то, учитывая (\ref{(3)}), получим
$$
H_n(x)=2\cdot{\rm Re}\left(\alpha_{\mu} e^{i {\mu}
x}\right)=\tau_{\mu}(x),\quad\hbox{где}
\quad\alpha_{\mu}=\frac{a_{\mu}}{2}-i\frac{b_{\mu}}{2}.
$$

\subsection{Случай, когда все $z_k$ различны}\label{s1.3}

В этом случае определитель системы (\ref{(4)}) --- определитель Вандермонда
$W(z_1,\ldots,z_n)\ne 0$, а ${Y}_k={Y}_k(\mu)$ находятся как
скалярные произведения (см., напр., \cite{ChuDanch},
\cite{DanchDod}): ${Y}_k=({\mathcal{L}}_k\cdot {\mathcal{S}})$,
где ${\mathcal{S}}= ({\sigma}_1,\ldots, {\sigma}_{n})$ и
\begin{equation}
{\mathcal{L}}_k=\frac{1}{\prod_{j\ne k}(z_k-z_j)}
\left((-1)^{n-1}\rho_{n-1}^{(k)}, \ldots ,
(-1)^{n-\mu}\rho_{n-\mu}^{(k)}, \ldots , -\rho_{1}^{(k)},\;
1\right),
  \label{(5)}
\end{equation}
где $\rho_m^{(k)}$ --- элементарные симметрические многочлены вида
$$
\rho_0=1,\quad \rho_m^{(k)}=\rho_m({z}_1,\ldots,{z}_{k-1},0,
{z}_{k+1},\ldots,{z}_n), \quad k=\overline{1,n},
$$
$$
\rho_m:=\rho_m({z}_1,\ldots,{z}_n)= \sum_{1 \le j_1<\ldots<j_m \le
n}{{z}_{j_1}\ldots {z}_{j_m}}, \quad m=\overline{1,n},
$$
Отсюда при условии $(\mathcal A)$ имеем
\begin{equation}
Y_k(\mu) =\frac{(-1)^{n-{\mu}}{\rho}^{(k)}_{n-{\mu}}}{\prod_{j\ne
k}(z_k-z_j)},\quad 1\le {\mu}\le n,\quad
X_k(1)=\frac{(-1)^{n-1}\prod_{j=1}^n z_j}{z_k^2 \cdot \prod_{l\ne
k}(z_k-z_j)}.
  \label{(6)}
\end{equation}

\subsection{Регулярное решение задачи (\ref{(2)}) при $\mu=1$ и $\mu=n$}\label{s1.4}

Приведем одно элементарное построение АФО.

\begin{example} \label{example1} При ${\mu}=n$ имеем регулярное решение:
\begin{equation}
\lambda_k=2\pi\frac{k-1}{n},\quad X_k=\frac{1}{n},\quad
k=\overline{1,n}.
 \label{(7)}
\end{equation}
При ${\mu}=1$ и нечетных $n\ge 3$ имеем регулярное решение $($с
точностью до знака $X_k$$)$:
\begin{equation}
\lambda_k=\frac{\pi(4k-n)}{n+2},\quad
X_k=\pm\frac{1}{2^{n-1}}{\prod_{s\ne
k}\left|\sin\frac{\lambda_s-\lambda_k}{2}\right|^{-1}},\quad
k=\overline{1,n}.
  \label{(8)}
\end{equation}
\end{example}

\begin{remark} \label{remark1} Несложно проверить, что в (\ref{(8)}) все
$z_k=\exp(-i\lambda_k)$ попарно различны и составляют множество
$$
\{\sqrt[n+2]{-1}\}\setminus \left\{e^{i\varphi},
e^{-i\varphi}\right\},\qquad \varphi:=\frac{n}{n+2}\pi.
$$
Ниже будут найдены все регулярные решения задачи (\ref{(2)}) при
$\mu=1$ независимо от четности $n$ (см. теорему \ref{theorem3}).
\end{remark}

\begin{remark} \label{remark1a}
Отметим, что в силу равенств $z_k^{n+2}=- 1$ нулевые степенные
суммы $\sum_k X_kz_k^{\beta}$ в (\ref{(4)}) порядка
$\beta\in\overline{2,n}$ остаются нулевыми при замене $\beta$ на
$\beta+{(n+2)}\,k$, $k\in{\mathbb N}$. Поэтому формула (\ref{(2)})
верна и для сходящихся тригонометрических рядов, в которых
отсутствуют гармоники с номерами $\pm 1+(n+2)k$, $(n+2)k$,
$k\in{\mathbb N}$.
\end{remark}

\begin{Proof} То, что (\ref{(7)}) действительно является решением
в случае $\mu=n$, легко проверяется непосредственно с учетом вида
$H_n(x)$ и равенств
$$
\sum_{k=1}^n \cos m\lambda_k=\sum_{k=1}^n \sin
m\lambda_k=\sum_{k=1}^n \sin n\lambda_k =0,\quad
m=\overline{1,n-1}.
$$

Перейдем к доказательству (\ref{(8)}). Из (\ref{(6)}) при $\mu=1$ имеем
\begin{equation}
(-1)^{\frac{n(n-1)}{2}}\prod_{j=1}^{n}X_j=
\frac{\prod_{s=1}^{n}z_s^{n-2}}{\prod_{1\le q<s\le
n}(z_q-z_s)^2}=\frac{1}{\prod_{s=1}^{n}z_s}\prod_{1\le q<s\le
n}\frac{z_q z_s}{(z_q-z_s)^2}.
   \label{(9)}
\end{equation}
Заметим, что последнее произведение вещественно, т.к.
$$
\frac{z_q z_s}{(z_q-z_s)^2}=-\frac{1}{4}{\rm cosec}^2
\frac{\lambda_q-\lambda_s}{2}, \quad z_q=e^{-i\lambda_q},\quad
z_s=e^{-i\lambda_s}.
$$
Поэтому для вещественности (\ref{(9)}) необходимо условие
\begin{equation}
{\rm Im}\,\prod_{s=1}^{n}z_s=0 \qquad\Leftrightarrow\qquad
\sum_{s=1}^{n}\lambda_s=0\,({\rm mod}\, \pi),\qquad
z_s=e^{-i\lambda_s}.
  \label{(10)}
\end{equation}
В силу (\ref{(10)}) числитель в (\ref{(6)}) должен быть вещественным. Значит, при
всех $k$ вещественным должен быть знаменатель:
\begin{equation}
v_k:=z_k^2 \cdot \prod_{s\ne k}(z_k-z_s),\quad {\rm
Arg}\,v_k=0\,({\rm mod}\, \pi).
   \label{(11)}
\end{equation}
Легко вычислить ${\rm Arg}\,
(z_k-z_s)=\frac{1}{2}(\pi+\lambda_s+\lambda_k)\;({\rm mod}\,
\pi)$, $|z_{k,s}|=1$. Поэтому аргумент в (\ref{(11)}) равен с точностью до
$\pi$ сумме
$$
2\lambda_k+\frac{1}{2}\sum_{s\ne
k}\lambda_s+\frac{n-1}{2}\lambda_k+\frac{n-1}{2}\pi=
\frac{1}{2}\sum_{s\ne
k}\lambda_s+\frac{n+3}{2}\lambda_k+\frac{n-1}{2}\pi.
$$
Теперь будем считать $n$ нечетным. Тогда с учетом предыдущего
имеем
\begin{equation}
{\rm Arg}\, v_k=\frac{1}{2}\Lambda+\frac{n+2}{2}\lambda_k=0
\;({\rm mod}\, \pi),\quad \Lambda:=\sum_{s=1}^n\lambda_s,
    \label{(12)}
\end{equation}
т.е. при всех $k=\overline{1,n}$ необходимо
$\Lambda+(n+2)\lambda_k=0\;({\rm mod}\, 2\pi)$. Рассмотрим
удовлетворяющую этому свойству систему
 \begin{equation}
\Lambda+(n+2)\lambda_k=4\pi k,\qquad k=\overline{1,n}.
  \label{(13)}
\end{equation}
Несложно найти явный вид ее решения. Действительно, после сложения
всех уравнений имеем: $(2n+2)\Lambda=2\pi n(n+1)$, откуда
$\Lambda=\pi n$. Таким образом, из (\ref{(13)}) имеем:
\begin{equation}
\pi n+(n+2)\lambda_k=4\pi k,\qquad
\lambda_k=\frac{\pi(4k-n)}{n+2},\qquad k=\overline{1,n}.
 \label{(14)}
\end{equation}
Подставив $\lambda_k$ в (\ref{(6)}), найдем $X_k=X_k(1)$. Покажем, что
найденные $\{\lambda_k\}$ и $\{X_k\}$ являются решением задачи
(\ref{(2)}). В самом деле, из (\ref{(13)}) следует (\ref{(12)}), откуда с учетом (\ref{(14)})
находим:
$$
{\rm Arg}\, v_k=\frac{1}{2}\Lambda+\frac{n+2}{2}\lambda_k=
\frac{\pi n}{2}+\frac{n+2}{2}\lambda_k=0\;({\rm mod}\,
\pi),\quad\hbox{т.е.} \quad v_k\in {\mathbb R}.
$$
Числители в (\ref{(6)}) также вещественны, точнее, они равны $-1$.
Действительно, из (\ref{(14)}) имеем
$$
\prod_{k=1}^n z_k=\exp\left(-i\sum_{k=1}^n
\lambda_k\right)=\exp\left(-i\sum_{k=1}^n
\frac{\pi}{n+2}(4k-n)\right)=e^{-i\pi n}=-1
$$
Для получения формулы (\ref{(8)}) остается вычислить
$$
|v_k|=\prod_{s\ne k}|e^{-i\lambda_s}-e^{-i\lambda_k}|=
2^{n-1}\prod_{s\ne k}
\left|\sin\frac{\lambda_j-\lambda_k}{2}\right|.
$$
\end{Proof}

\newpage
\section{Сведение к задаче дискретных моментов}\label{s2}
\subsection{Метод Прони}\label{s2.1}
Дополним (\ref{(4)}) до системы
\begin{equation}
\left\{
\begin {array}{lllllllll}
 X_1 z_1^{-n+1}&+&X_2 z_2^{-n+1}&+&\ldots&+& X_{n} z_{n}^{-n+1} &=
 &{\sigma}_{1-n}\\
 \ldots &&&& \ldots &&&&\ldots\\
X_1 z_1^{-1}&+&X_2 z_2^{-1}&+&\ldots&+& X_{n} z_{n}^{-1} &= &{\sigma}_{-1}\\
X_1 &+&X_2 &+&\ldots&+& X_{n} &= &{\sigma}_0={\omega}\\
 X_1 z_1&+&X_2 z_2&+&\ldots&+& X_{n} z_{n} &= &{\sigma}_1\\
\ldots &&&& \ldots &&&&\ldots\\
 X_1 z_1^n&+&X_2 z_2^n&+&\ldots&+& X_{n} z_{n}^n &= &{\sigma}_n
\end {array}
\right.
 \label{(15)}
\end{equation}
где $z_k\ne 0$ и, возможно, некоторые $X_k=0$. Для разрешимости
задач (\ref{(2)}), (\ref{(4)}) при условиях $({\mathcal A})$
необходима и достаточна разрешимость задачи (\ref{(15)}) (с
каким-либо вещественным $\omega$) при условиях

\smallskip
 {$({\mathcal B}')$}:
 \textit{
 ${\sigma}_{0}=\omega$,
 ${\sigma}_{\mu}={\sigma}_{-{\mu}}=1$,
 ${\sigma}_k=0$ при $k\ne 0,\pm {\mu}$;\;
все $z_k$ различны и $|z_k|=1$. }

\smallskip
Необходимость очевидна, поскольку добавленные уравнения получаются
комплексным сопряжением уравнений (\ref{(4)}). Достаточность
следует из того, что недостающее для $({\mathcal A})$ условие
$X_k\in {\mathbb R}$ является следствием условий $({\mathcal
B}')$. Справедлива

\begin{lemma} \label{lemma1}Для совместной системы $(\ref{(15)})$ условие $X_k\in
{\mathbb R}$ следует из условий $({\mathcal B}')$.\end{lemma}

\begin{Proof} Пусть для определенности $X_k\ne 0$ при
$k=\overline{1,m}$, $m\le n$, а остальные $X_k$ равны нулю.
Заметим, что операция сопряжения к уравнениям с номерами
$\overline{n-m+1,n}$ даст уравнения с номерами
$\overline{n,n+m-1}$ с заменой $X_k$ на $\overline{X_k}$. Другими
словами, уравнения с номерами $\overline{n,n+m-1}$ верны и для
$X_k$, и для $\overline{X_k}$. Вычитая эти уравнения с одинаковыми
номерами, получим систему относительно $X_k-\overline{X_k}$ с
нулевыми правыми частями и определителем Вандермонда
$W(z_1,\ldots,z_m)$ с попарно различными $z_k$. Она имеет нулевое
решение.\end{Proof}

В ряде выкладок мы будем изучать (\ref{(15)}) без каких-либо
предварительных условий на $\{X_k\}$, $\{z_k\}$, $\{{\sigma}_k\}$,
$\omega$. Это классическая система дискретных моментов
относительно неизвестных $\{X_k\}$, $\{z_k\}$ при прочих известных
параметрах (см. например, \cite{Prony}, \cite{Sylvester},
\cite{Kung1}). Вместе с (\ref{(15)}) будем рассматривать систему моментов
\begin{equation}
\omega_l=S_l,\quad \hbox{где}\quad \omega_l:=\sum_{k=1}^n Z_k
z_k^l,\qquad l=\overline{0,2n-1},
 \label{(16)}
\end{equation}
относительно $Z_k$, $z_k$ при известных $S_l={\sigma}_{1-n+l}$,
$(S_{n-1}=\omega)$. В случае, когда $z_k\ne 0$,
$k=\overline{1,n}$, система (\ref{(16)}) получается из (\ref{(15)}) заменой
$Z_k=X_k z_k^{-n+1}$, так что в этом случае системы равносильны.
Отметим, что в задаче (\ref{(16)}) условие $z_k\ne 0$ необязательно.

\subsection{Регулярные системы моментов}\label{s3.2}

Следуя \cite{Lyubich}, совместную систему (\ref{(16)}) и ее решение будем
называть \textit{регулярными}, если выполнены условия

\smallskip
\textbf {$({\mathcal B}'')$}:   \textit {все $z_k$ попарно
различны $($не обязательно отличны от нуля$)$, а все $Z_k$ отличны
от нуля.}

\smallskip
 Положим $({\mathcal B})=({\mathcal B}')\cup ({\mathcal
B}'')$ --- совокупность условий $({\mathcal B}')$ и $({\mathcal
B}'')$, т.е.

\smallskip
\textbf {$({\mathcal B})$}:   \textit
{${\sigma}_{0}=\omega\in{\mathbb R}$,
 ${\sigma}_{\mu}={\sigma}_{-{\mu}}=1$,
${\sigma}_k=0$ при $k\ne 0,\pm {\mu}$; все $z_k$ попарно различны
и $|z_k|=1$, все $Z_k$ отличны от нуля.}

\smallskip
Для краткости совместные системы с указанными свойствами
$({\mathcal B}')$, $({\mathcal B}'')$, $({\mathcal B})$ и их
решения будем называть $({\mathcal B}')$-разрешимыми, $({\mathcal
B}'')$-регулярными, $({\mathcal B})$-регулярными соответственно.

Для нахождения решения $({\mathcal B}'')$-регулярных систем
(\ref{(16)}) обычно пользуются классическим методом Прони
\cite{Prony}. В нашем случае он состоит в следующем. Рассмотрим
два определителя
$$
d_1=\left|
\begin{array}{ccccc}
1 & 0 & 0 & \ldots & 0\\
0 & Z_1 & Z_2 & \ldots & Z_n\\
0 & Z_1{z}_1 & Z_2{z}_2 & \ldots & Z_n{z}_n\\
\ldots & \ldots & \ldots & \ldots & \ldots\\
0 & Z_1{z}_1^{n-1} & Z_2{z}_2^{n-1} & \ldots & Z_n{z}_n^{n-1}\\
\end{array}
\right|,
 \quad
d_2(z)=\left| \begin {array}{ccccc}
  1&z&{z}^{2}&\ldots&{z}^{n}
\\\noalign{\medskip}
  1&z_1&z_1^{2}&\ldots&z_1^{n}
\\\noalign{\medskip}
   \ldots&\ldots&\ldots&\ldots&\ldots
   \\\noalign{\medskip}
  1&z_n&z_n^{2}&\ldots&z_n^{n}
\end {array} \right|.
$$
Тогда в силу равенств (\ref{(16)}) с учетом замены $S_l={\sigma}_{1-n+l}$,
$S_{n-1}={\sigma}_0=\omega$, имеем
\begin{equation}
d_1\cdot d_2(z)=\left| \begin {array}{cccccccc}
1&z&z^2&\ldots&\ldots&z^{n-2}& z^{n-1}&z^n
\\\noalign{\medskip}{\sigma}_{1-n}&{\sigma}_{2-n}&\ldots&\ldots&{\sigma}_
{-2}&{\sigma}_{-1}&{\omega}&{\sigma}_1
\\\noalign{\medskip}{\sigma}_{2-n}&\ldots&\ldots&{\sigma}_{-2}
&{\sigma}_{-1}&{\omega}&{\sigma}_1&{\sigma}_2
\\\noalign{\medskip}\ldots&\ldots
&\ldots&\ldots&\ldots&\ldots&\ldots&\ldots
\\\noalign{\medskip}{\sigma}_{-1}&{\omega}&{\sigma}_1&{\sigma}_2&\ldots&\ldots&{\sigma}_{n-2}
&{\sigma}_{n-1}
\\\noalign{\medskip}{\omega}&{\sigma}_1&{\sigma}_2
&\ldots&\ldots&{\sigma}_{n-2}&{\sigma}_{n-1}&
{\sigma}_n\end {array} \right|
 \label{(17)}
\end{equation}
Из условия $({\mathcal B}'')$-регулярности системы (\ref{(16)})
следует, что определитель $d_1$ отличен от нуля, а $d_2(z)$
обращается в нуль только при $z=z_k$ (определители Вандермонда).
Следовательно, искомые числа $z_k$ являются простыми корнями
многочлена $G_n(\{\sigma_k\};\,z):=d_1\cdot d_2(z)$. Будем
называть
$$
G_n(\{\sigma_k\};\,z)=g_0+g_1 z+\ldots+ g_n z^n,\qquad
g_k=g_k(\omega),
$$
{\it производящим многочленом}. Решив уравнение
$G_n(\{\sigma_k\};\,z)=0$ относительно $z$, найдем $z_k$, а затем,
используя первые $n$ уравнений системы (\ref{(16)}) (с отличным от нуля
определителем Вандермонда), по формулам (\ref{(6)}) определим $Z_k$.

Вообще говоря, представление (\ref{(17)}) для
$G_n(\{\sigma_k\};\,z)$ будем применять без каких-либо
предварительных условий на $\sigma_k$. В частности, степень $G_n$
может быть меньше $n$, возможно даже, что
$G_n(\{\sigma_k\};\,z)\equiv 0$. При выполнении условий
$({\mathcal B}')$ на $\sigma_k$ вместо $G_n(\{\sigma_k\};\,z)$ для
определенности будем использовать обозначение
$G_n(\omega,\mu;\,z)$. В этом случае в (\ref{(17)}) под первой
строкой получается трехдиагональная матрица с элементами
$$
A_{k,n+2-k}=\omega,\quad A_{k,n+2-k\pm \mu}=1,\quad k\ge 2
$$
(в этих равенствах соответственно $k=\overline{2,n+1}$,
$k=\overline{\mu+1,n+1}$, $k=\overline{2,n+1-\mu}$), остальные
элементы равны нулю. Например,
$$
G_4(\omega,1;\,z) =
 \left|\begin {array}{ccccc} 1&z&{z}^{2}&{z}^{3}&{z}^{4}
\\\noalign{\medskip}0&0&1&{\omega}&1\\\noalign{\medskip}0&1&{\omega}&1&0
\\\noalign{\medskip}1&{\omega}&1&0&0\\\noalign{\medskip}{\omega}&1&0&0&0\end {array}
 \right|,\quad
G_4(\omega,2;\,z) =
 \left|\begin {array}{ccccc} 1&z&{z}^{2}&{z}^{3}&{z}^{4}
\\\noalign{\medskip}0&1&0&{\omega}&0\\\noalign{\medskip}1&0&{\omega}&0&1
\\\noalign{\medskip}0&{\omega}&0&1&0\\\noalign{\medskip}{\omega}&0&1&0&0\end {array}
 \right|.
$$
Приведем известный критерий $({\mathcal B}'')$-регулярности систем
(\ref{(16)}).

\begin{thS-L}[{\cite{Sylvester},\cite{Lyubich}}] \label{thS-L}
Для  $({\mathcal B}'')$-регулярности системы $(\ref{(16)})$
необходимо и достаточно, чтобы степень производящего многочлена
$G_n(\{\sigma_k\};\,z)$ была равна $n$, а все его корни
$z_1,\ldots,z_n$ были попарно различны. Они являются решением
системы $(\ref{(16)})$. При этом регулярная система имеет
единственное решение.
\end{thS-L}

В связи с этой теоремой отметим, что в части необходимости она
следует из приведенного метода Прони. Достаточность гарантирует,
что если корни производящего многочлена попарно различны, то
система (\ref{(16)}) совместна и притом с отличными от нуля $Z_k$.
Связано это с тем, что если какой-либо $Z_k=0$ в совметной системе
(\ref{(16)}), то определитель (\ref{(17)}) равен тождественно нулю
(подробнее, см. в п.~\ref{s4.1}).

\begin{lemma} \label{lemma2} В случае $({\mathcal B}'')$-регулярной
разрешимости задачи $(\ref{(16)})$ производящий многочлен, построенный по
формуле $(\ref{(17)})$, имеет вид
\begin{equation}
G_n(\{\sigma_k\};\,z)=(-1)^n\prod_{k=1}^n Z_k\cdot \prod_{1\le
k<j\le n} (z_k-z_j)^2 \cdot\prod_{k=1}^n (z-z_k)
  \label{(18)}
\end{equation}
\end{lemma}

\begin{Proof} То, что в (\ref{(18)}) присутствует последнее
произведение следует из теоремы Сильвестра-Любича. Хорошо известно
(см., например, \cite{Boley}), что матрица
$$
A:=\left(\begin {array}{cccc}
 \omega_0&\omega_1&\ldots&\omega_{n-1}\\
\noalign{\medskip}
 \omega_1&\omega_2&\ldots&\omega_{n-2}\\
\noalign{\medskip}
 \ldots&\ldots&\ldots&\ldots\\
  \omega_{n-1}&\omega_n&\ldots&\omega_{2n}\\
\end{array} \right),\qquad \omega_l:=\sum_{k=1}^n Z_k
z_k^l,
$$
получающаяся из матрицы (\ref{(17)}) вычеркиванием первой строки и
последнего столбца представляется в виде произведения матрицы
Вандермода, диагональной матрицы и транспонированной матрицы
Вандермода:
$$
A=W\cdot{\rm diag}(Z_1,\ldots,Z_n)\cdot W^{T},\qquad W:= \left(
\begin {array}{ccc}
1&\ldots&1\\
\noalign{\medskip}z_{{1}}&\ldots&z_{{n}}\\
\noalign{\medskip}\ldots&\ldots&\ldots\\
\noalign{\medskip}{z_{1}}^{n-1}&\ldots&z_n^{n-1}
\end{array} \right)
$$
То есть алгебраическое дополнение определителя (\ref{(17)}) к $z^n$
(старший коэффициент многочлена $G_n(\{\sigma_k\};\,z)$) равно
$(-1)^n (\det W)^2 \prod_{k=1}^n Z_k$.
\end{Proof}

\section{О некоторых свойствах  производящего многочлена}\label{s3}
\subsection{}\label{s3.1}

Докажем простое утверждение.

\begin{lemma} \label{lemma3} Пусть $G_n(\{\sigma_k\};\,z)$ --- вещественный
многочлен степени $n$, корни которого лежат на единичной
окружности $($симметрично относительно действительной оси$)$.
Тогда коэффициенты этого многочлена симметричны или
антисимметричны, т.е. с одним из знаков $\pm$ имеем:
\begin{equation}
g_k=\pm g_{n-k},\qquad g_k=g_k(\omega),\qquad k=\overline{0,n}.
 \label{(19)}
\end{equation}
\end{lemma}
\begin{Proof} В самом деле, поскольку все $|z_k|=1$, то  $g_n=\pm g_0$ (теорема
Виета для вещественных коэффициентов). Далее, при $|z|=1$ имеем
$$
\overline{G_n(\{\sigma_k\};\,z)}=G_n(\{\sigma_k\};\,\overline{z})=
G_n(\{\sigma_k\};\,z^{-1})=z^{-n} G(z),
$$
где $G$ --- многочлен с коэффициентами $\tilde g_k=g_{n-k}$,
расположенными в обратном порядке. Следовательно, многочлены $G$ и
$G_n$ имеют одинаковые корни, а их старшие коэффициенты могут
отличаться только знаком. Отсюда получается (\ref{(19)})
\end{Proof}

\begin{example} \label{example2} В некоторых случаях условие симметрии
коэффициентов сразу приводит к решению задачи (\ref{(2)}). Пусть
$n=4$, $\mu=2$. Тогда
$$
G_4(\omega,2;\,z)=
(1-\,{\omega}^{2})\cdot((1-\,{\omega}^{2})z^4+\omega z^2-1).
$$
Из условия $g_0(\omega)=g_4(\omega)$ находим $\omega=\sqrt{2}$.
Подставив это значение в $G_4$, найдем корни, а затем по формулам
(\ref{(6)}) --- амплитуды:
$$
z_k=e^{-i\lambda_k}=\pm\frac{1}{2}\,\sqrt {2-\sqrt
{2}}\pm\frac{i}{2}\,\sqrt {2+\sqrt {2}},\quad
X_1=\ldots=X_4=\frac{\sqrt{2}}{4}.
$$
Отсюда находим регулярное решение  $\{\lambda_k\}=\{\pm
3\pi/8,\,\pm 5\pi/8\}$ задачи (\ref{(2)}) и АФО $H_4$. Аналогично
при $n=6$, $\mu=2$ условие $g_0(\omega)=g_6(\omega)$ приводит к
значению $\omega=(\sqrt{5}-1)/2$ и нужному АФО $H_6$.
\end{example}

\subsection{Рекуррентные формулы для коэффициентов многочлена $G_n(\omega,\mu;\,\cdot)$}\label{s3.2}
Заметим, что при условиях $({\mathcal B}')$ на $\sigma_k$ строки
$T_j$ определителя (\ref{(17)}) с номерами $j=2,\ldots, n+1$ имеют
вид
\begin{equation}
T_j=({\sigma}_{n-j+1},\ldots,{\sigma}_{\mu},
\ldots,{\sigma}_1,{\omega},{\sigma}_1,\ldots,{\sigma}_{\mu},
\ldots,{\sigma}_{j-1}),
 \label{(20)}
\end{equation}
где ${\sigma}_{\mu}=1$, все ${\sigma}_k=0$ с номерами $k\ne
0,{\mu}$, а ${\sigma}_0=\omega$
--- вещественный вспомогательный параметр (некоторые ${\sigma}_{k}$ могут
выходить за пределы строк и тогда эти элементы отбрасываются).

Запишем коэффициенты многочлена $G_n(\omega,\mu;\,\cdot)$ в виде
строки
$$
{\mathcal G}:=(g_0,\,g_1,\ldots,g_n).
$$
Хорошо известна следующая обобщенная формула Ньютона: скалярное
произведение $(T_j\cdot {\mathcal G})=0$ (здесь она получается как
сумма произведений алгебраических дополнений к элементам первой
строки определителя (\ref{(17)}) на соответствующие элементы
другой строки).

Далее, считая для удобства, что $g_k=0$ с номерами $k<0$ и $k>n$
перепишем это условие в развернутом виде:
$$
g_{k-{\mu}}+g_{k}\cdot {\omega}+g_{k+{\mu}}=0,\qquad
k=\overline{0,n-1},
$$
откуда видим, что эта система распадается на ${\mu}$ независимых
подсистем
\begin{equation}
g_{q+{\mu} k}+{\omega} g_{q+{\mu}(k-1)}+g_{q+{\mu}(k-2)}=0,\qquad
q=\overline{0,\mu-1},\quad k\in{\mathbb N}.
 \label{(21)}
\end{equation}
Введем рекуррентно вспомогательные многочлены $r_s(\omega)$:
\begin{equation}
r_{-1}=0,\quad r_0=1,\quad r_{k}(\omega)+{\omega}
r_{k-1}(\omega)+r_{k-2}(\omega)=0,\quad \quad k\in{\mathbb N}.
 \label{(22)}
\end{equation}
С учетом равенств (\ref{(21)}), (\ref{(22)}) получается

\begin{lemma} \label{lemma4} Коэффициенты производящего многочлена
$G_n(\omega,\mu;\,\cdot)$ при условиях $({\mathcal B}')$ на
$\sigma_k$ линейно выражаются через $g_0,\ldots,g_{{\mu}-1}$ и
имеют вид
\begin{equation}
g_{q+{\mu}k}=g_q r_{k},\qquad q=\overline{0,\mu-1}, \quad
k=0,1,\ldots \quad (q+{\mu}k\le n).
  \label{(23)}
\end{equation}
\end{lemma}

Хорошо известно \cite{SUETIN}, что равенствами (\ref{(22)}) при
$\omega=-2x$ определяются многочлены Чебышева $U_{k}$ второго
рода. Точнее,
\begin{equation}
U_{k}(x):=\frac{\sin((k+1)\arccos x)}{\sqrt{1-x^2}}= r_k(-2 x);
\qquad r_k(\omega)=U_{k}\left(-\frac{\omega}{2}\right).
 \label{(24)}
\end{equation}
Явный вид таких многочленов \cite{SUETIN}:
\begin{equation}
r_{k}({\omega})=(-1)^k \sum_{j=0}^{k/2}(-1)^j C_{k-j}^j
\omega^{k-2j}.
  \label{(25)}
\end{equation}
Отсюда видно, что
\begin{equation}
r_{k}({\omega})=(-1)^kr_{k}({-\omega}),\qquad
r_{k}({\omega})=(-1)^kU_{k}\left({\omega}/{2}\right).
   \label{(26)}
\end{equation}
 Кроме того, из (\ref{(24)})
следует, что многочлены $r_k$ образуют ортонормированную систему в
$L_2([-2,2])$ с весом $(2\pi)^{-1}\sqrt{4-\omega^2}$.

Приведем одно необходимое условие регулярной разрешимости задачи
(\ref{(2)}) в терминах $r_k$. Пусть
$$
s=\left[{n}/{{\mu}}\right],\qquad n-{\mu}<{\mu}  s\le n
$$
и предположим, что число $n-{\mu} s$ четное. При $q=(n-{\mu} s)/2$
соответствующая последовательнось (\ref{(23)}) содержит два коэффициента с
симметричными номерами $g_q$ и $g_{n-q}$. Отсюда и из (\ref{(19)}), (\ref{(23)}),
(\ref{(24)}) получается

 \begin{theorem} \label{theorem1} В случае четного $n-{\mu} s$ для $({\mathcal
B})$-регулярной разрешимости системы $(\ref{(15)})$ необходимо одно из
равенств
\begin{equation}
U_{s}\left({\omega}/{2}\right)=\pm 1.
    \label{(27)}
\end{equation}
\end{theorem}

 \begin{example} \label{example3} Пусть $n=8$, $\mu=2$. В этом случае $ s=4$.
Выпишем уравнение (\ref{(27)}) со знаком "+" и его корни:
$$
{\omega}^{4}-3\,{\omega}^{2}=0;\qquad \{0,\,0,\,\pm\sqrt{3}\}.
$$
Непосредственной подстановкой каждого из этих корней в
$G_8(\omega,2;\,z)$ убеждаемся, что условия $({\mathcal B})$ по
отношению к корням многочлена $G_8$ выполняются. Например, при
подстановке $\omega=0$ получаем $G_8(0,2;\,z)={z}^{8}-{z}^{4}+1$ с
корнями $\{z_k\}=\{e^{-i\lambda_k}\}=\sqrt[4]{\exp(\pm i\,
\pi/3)}$. Значит, задача (\ref{(2)}) во всех трех случаях имеет
регулярные вещественные решения $\{X_k\}$ (см.
лемму~\ref{lemma1}), которые находятся из (\ref{(6)}).
\end{example}

В дальнейшем важную роль играет многочлен
\begin{equation}
{\mathcal R}_{s-1}(\omega;t):=r_0(\omega)+r_{1}(\omega) t+\ldots+
r_{s-1}(\omega)t^{s-1},\quad r_0(\omega)\equiv 1,\quad s\in
{\mathbb N}.
   \label{(28)}
\end{equation}

\subsection{О корнях многочлена ${\mathcal R}_{s-1}$}\label{s3.3}

Справедлива
 \begin{lemma} \label{lemma5} Все корни $t_k$ многочлена ${\mathcal R}_{s-1}$,
$s\ge 2$, различны и лежат на единичной окружности в том и только
том случае, когда $\omega$ является корнем многочлена
$r_{s}({\omega})=(-1)^sU_{s}\left({\omega}/{2}\right)$, т.е.
\begin{equation}
\omega\in\Omega_s:=\left\{2\cos
\varphi_{\alpha}:\quad\varphi_{\alpha}
 =\frac{\pi {\alpha}}{s+1},\quad {\alpha}=\overline{1,s}\right\}.
   \label{(29)}
\end{equation}
При  $\omega=-2\cos \varphi_{\alpha}=2\cos \varphi_{s-\alpha+1}$
корни $t_k$ составляют множество
\begin{equation}
\{t_k\}=\{\sqrt[s+1]{(-1)^{\alpha}}\}\setminus
\{e^{i\varphi_{\alpha}},e^{-i\varphi_{\alpha}}\}.
   \label{(30)}
\end{equation}
\end{lemma}

\begin{Proof} При $x=-\omega/2$, $\varphi=\arccos x$ из
(\ref{(24)}) имеем
$$
{\sqrt{1-x^2}}{\mathcal
R}_{s-1}(\omega;t)={\sqrt{1-x^2}}\sum_{k=0}^{s-1}t^sU_{k}\left(x\right)=
\sum_{k=0}^{s-1}t^s \sin((k+1)\varphi)=
$$
$$
\frac{1}{2i}\sum_{k=0}^{s-1}t^s \left(e^{i (k+1)\varphi}-e^{-i
(k+1)\varphi}\right )= \frac{1}{2i}\left(e^{i \varphi}\frac{1-t^s
e^{i s\varphi}}{1-t e^{i\varphi}}-e^{-i \varphi}\frac{1-t^s e^{-i
s\varphi}}{1-t e^{-i\varphi}}\right)=
$$
\begin{equation}
=\frac{t^{s+1} \sin(s\varphi)-t^s
\sin((s+1)\varphi)+\sin(\varphi)}{t^2-2 t\cos \varphi +1}=:P(t).
 \label{(31)}
\end{equation}
Для того, чтобы многочлен в числителе последней дроби имел
различные корни $t_k$ с $|t_k|=1$ по лемме 3 необходимо (а в
данном случае и достаточно) чтобы
$$
\sin(s+1)\varphi=0 \hbox{ и } \sin s\varphi=\sin\varphi \hbox{
~~~или~~~ } \sin(s+1)\varphi=0 \hbox{ и } \sin
s\varphi=-\sin\varphi.
$$
Решив эти две системы уравнений, получим (первой и второй системам
соответствуют следующие первое и второе равенства):
\begin{equation}
\cos \frac{s+1}{2}\varphi=0\quad\hbox {или}\quad \sin
\frac{s+1}{2}\varphi=0,\qquad P(t)=A \cdot \frac{t^{s+1} \pm
1}{t^2-2 t\cos \varphi +1}.
    \label{(32)}
\end{equation}
Здесь $A=\sin(s\varphi)\ne 0$. Отсюда получается (\ref{(29)}). Исключив из
множества корней числителя в $P$ нули знаменателя, получим (\ref{(30)}).
\end{Proof}

Отметим, что из ортогональности многочленов $r_k$, в частности,
следует, что (вещественные) корни многочленов $r_{s-1}$ и $r_{s}$
все различны \cite{SUETIN}, так что ${\rm deg}\,{\mathcal
R}_{s-1}(\omega;t)=s-1$ при любых $\omega$, лежащих в достаточно
малой окрестности множества $\Omega_s$. Из сказанного получается

 \begin{lemma} \label{lemma6} Пусть $\omega_0\in \Omega_s$. Тогда при
достаточно малых $|\varepsilon|$ многочлен ${\mathcal
R}_{s-1}(\omega_0+\varepsilon;t)$ имеет ровно $s-1$ различных
корней и при $|\varepsilon|\to 0$ они стремятся к корням
многочлена ${\mathcal R}_{s-1}(\omega_0;t)$.
\end{lemma}

\subsection{О факторизации производящего многочлена при условиях
$({\mathcal B}')$ на $\sigma_k$}\label{s3.4}
Рассмотрим случай
$n={\mu}s-1$, $s\in{\mathbb N}$, $s\ge 2$. Из (\ref{(23)}) получается
\begin{equation}
G_{{\mu}s-1}(\omega,\mu;\,z)=(g_0+g_1 z+\ldots+g_{{\mu}-1}
z^{{\mu}-1})\cdot {\mathcal R}_{s-1}(\omega;\,z^{\mu}).
    \label{(33)}
\end{equation}
Покажем, что первый сомножитель может содержать лишь одно отличное
от нуля слагаемое $g_{{\mu}-1} z^{{\mu}-1}$. Справедлива

 \begin{lemma} \label{lemma7} Коэффициенты $g_{m+k\mu}$, $k\in {\mathbb
N}$, производящего многочлена $G_n(\omega,\mu;\,z)$ в случае
$n={\mu}s-1$ равны нулю при $m=\overline{0,\mu-2}$.
 \end{lemma}

\begin{Proof} Рассмотрим матрицу, получающуюся удалением
первой строки из матрицы в (\ref{(17)}). Покажем, что ее столбцы
$a_{k\mu}$ с номерами $k\mu$, $k=\overline{1,s}$, линейно
зависимы. Отсюда и будет следовать утверждение леммы.
Действительно, запишем равную нулю линейную комбинацию $\sum
x_{k\mu} a_{k\mu}=0$ в развернутом виде (ср. с (\ref{(21)}), где $q=\mu$):
$$
  x_{(k+1)\mu}+{\omega}
x_{k\mu}+x_{(k-1){\mu}}=0,\quad x_0=0,\; x_{\mu}=1,\quad
k=\overline{1,s-1}.
$$
Таким образом нетривиальная нулевая комбинация получается с
коэффициентами $x_{k\mu}$, определяемыми рекуррентно.
\end{Proof}

Отсюда и из леммы~\ref{lemma4} получается

\begin{theorem} \label{theorem2} Пусть $n={\mu}s-1$, $\mu\ge 1$, $s\ge 2$. При
условиях $({\mathcal B}')$ на $\sigma_k$ коэффициенты
производящего многочлена $G_n(\omega,\mu;\,z)$ имеют вид
\begin{equation}
g_{{\mu}k-1}=g_{\mu-1} r_{k-1},\qquad k=\overline{1,s},
  \label{(34)}
\end{equation}
остальные коэффициенты равны нулю и $(\ref{(33)})$ принимает вид
 \begin{equation}
  G_{n}(\omega,\mu;\,z)=g_{{\mu}-1}(\omega) z^{{\mu}-1} \cdot {\mathcal
R}_{s-1}(\omega;z^{\mu}).
   \label{(35)}
\end{equation}
\end{theorem}

Отметим, что $g_{\mu-1}(\omega)$ --- многочлен степени $n$, не
равный тождественно нулю. Действительно, при $\omega\to\infty$
старшее слагаемое определителя $G_{n}(\omega,\mu;z)$ имеет вид
$\omega^{n}(1+o(1))z^{n}$. Из (\ref{(35)}), в частности, следует,
что при $\mu\ge 2$ все решения задач (\ref{(2)}) и (\ref{(4)}) не
являются регулярными (если они существуют).

Приведем следствие, дополняющее пример~\ref{example1}.

 \begin{theorem} \label{theorem3} При $n\ge 2$ и $\mu=1$ для регулярной
разрешимости задачи $(\ref{(15)})$ $($а значит, и $(\ref{(2)}))$
необходимо и достаточно условие $\omega\in \Omega_s$ при $s=n+1$.
При $\omega=-2\cos \varphi_{\alpha}$ решение $z_k$
$(k=\overline{1,n})$ системы $(\ref{(4)})$ составляет множество
$(\ref{(30)})$ $($с $s=n+1)$ и
$$
X_k=-\frac{{\rm
sgn}(\sin(n+1)\varphi_{\alpha})}{n+2}\cdot\frac{z_k^2-2z_k\cos
\varphi_{\alpha}+1}{z_k^{n+3}},\quad
\varphi_{\alpha}=\frac{\pi}{n+2}{\alpha},\quad
\alpha=\overline{1,n+1}.
$$
\end{theorem}
\begin{Proof}
Действительно, при $\mu=1$ и условиях $({\mathcal B}')$ найдем
алгебраическое дополнение к $z^0$ в (\ref{(17)}). Его разложение
по последней строке, содержащей единственный ненулевой элемент
$\sigma_1=1$ в левом нижнем углу, дает определитель той же
структуры, имеющий на единицу меньший размер. Далее все
повторяется и по индуции находим $g_0=(-1)^{n+1}\ne 0$. Значит, из
(\ref{(35)}) имеем $G_{n}(\omega,1;\,z)=(-1)^{n+1} {\mathcal
R}_{n}(\omega;z)$, так что все корни производящего многочлена при
$\omega\in \Omega_{n+1}$ различны и лежат на единичной окружности.
Поэтому первая часть следствия получается из теоремы
Сильвестра-Любича.

Значения $X_k$ получаются из (\ref{(6)}), если учесть, что
знаменатель выражения (\ref{(6)}) равен $(-1)^{n-1}z_k^2{\mathcal
P}'(z_k)$, где ${\mathcal P}=P/A$ --- многочлен, определенный в
(\ref{(32)}).
\end{Proof}

\section{Разрешимая задача в нерегулярном случае}\label{s4}
\subsection{Регуляризация. Первый способ}\label{s4.1}

В ряде случаев система (\ref{(15)}) может быть разрешимой, но не
$({\mathcal B}'')$-регулярной (т.е. некоторые $X_k$ равны нулю и
все $z_k\ne 0$). В таких случаях, сделав замену переменных как в
(\ref{(16)}), ее можно переписать в виде переопределенной
совместной системы моментов c $m<n$:
\begin{equation}
\sum_{k=1}^m Z_k z_k^l=S_l,\qquad l=\overline{0,2n-1},\quad
Z_k=X_k z_k^{-n+1},\quad S_l={\sigma}_{1-n+l},
 \label{(36)}
\end{equation}
где все $Z_k$ отличны от нуля, а $z_k$ попарно различны (в
противном случае чисто неизвестных в  (\ref{(36)}) будет меньше
$m$). Отметим, что решение совместной системы (\ref{(36)})
единственно. Действительно, достаточно взять подсистему уравнений
с номерами $l=\overline{0,2m-1}$ и воспользоваться теоремой
Сильвестра-Любича.

Итак, пусть система (\ref{(36)}) совместна при $m<n$. В этом случае
производящий многочлен $G_n(\{\sigma_k\};z)$, построенный по
формуле (\ref{(17)}) (без каких-либо условий на $\sigma_k$), тождественно
равен нулю. Это следует из леммы~\ref{lemma2}. Действительно, достаточно
дополнить систему (\ref{(36)}) до регулярной добавлением к ее обеим частям
слагаемых с малыми $Z_k=\varepsilon$:
$$
\sum_{k=1}^{m} Z_k z_k^l+
\varepsilon\sum_{k=m+1}^{n}{z}_k^{l}\equiv
\xi_{l}(\varepsilon);\qquad
\xi_{l}(\varepsilon)=S_l+O(\varepsilon),
   \eqno{(37)}
$$
где все $z_k$ различны и фиксированы. Тогда при $\varepsilon\to 0$
коэффициенты соответствующго производящего многочлена в силу (\ref{(18)})
стремятся к нулю, они же в силу формулы (\ref{(17)}) при
$\xi_{l}(\varepsilon)\to S_l$ стремятся и к коэффициентам
многочлена $G_n$. Отсюда  получается нужное.

Далее, при малых $\varepsilon>0$ и $\varepsilon_k\in
(0,\varepsilon)$ дополним систему (\ref{(36)}) до регулярной
добавлением к ее обеим частям степенных слагаемых:
    \begin{equation}
\sum_{k=1}^{m} Z_k z_k^l+
\sum_{k=1}^{n-m}{\varepsilon}_k^{l}\equiv \xi_{l};\qquad
\xi_{l}=\xi_{l}(\{\varepsilon_k\})=S_l+O(\varepsilon),
   \label{(38)}
    \end{equation}
где все $z_k$ и $\varepsilon_k$ попарно различны и отличны от
нуля. Для этой системы также построим производящий многочлен (см.
(\ref{(18)}))
\begin{equation}
G_n(\{\xi_k\};z)=(-1)^nK(\{\varepsilon_k\})\cdot \prod_{k=1}^m
Z_k\cdot
 \prod_{k<j} (z_k-z_j)^2
\cdot\prod_{k=1}^m (z-z_k)\cdot\prod_{k=1}^{n-m}
(z-\varepsilon_k),
  \label{(39)}
\end{equation}
где
\begin{equation}
K(\{\varepsilon_k\})=\prod_{k<j}
(\varepsilon_k-\varepsilon_j)^2\cdot
\prod_{k=1}^{m}\prod_{j=1}^{n-m} (z_k-\varepsilon_j)^2.
  \label{(40)}
\end{equation}
Положим
$$
G_n^*(z)=\lim_{\varepsilon\to 0} \frac{G_n(\{\xi_k\};z)}{
\prod_{k<j} (\varepsilon_k-\varepsilon_j)^2}=
$$
\begin{equation}
=(-1)^n\prod_{k=1}^m Z_k\cdot\prod_{k<j}
(z_k-z_j)^2\cdot\prod_{k=1}^{m}z_k^{2(n-m)} \cdot z^{n-m}
\prod_{k=1}^m (z-z_k).
   \label{(41)}
\end{equation}
Будем называть $G_n^*(z)$ {\it нерегулярным производящим многочленом}.
Итак, справедлива

 \begin{theorem} \label{theorem4} Если переопределенная система $(36)$
совместна $($при отличных от нуля $z_k)$, то необходимо
$G_n(\{\sigma_k\};z)\equiv 0$, кроме того, $z_k$ являются
ненулевыми корнями нерегулярного производящего многочлена
$G_n^*(z)$.
\end{theorem}

Отсюда получается следующий

\smallskip
{\bf Алгоритм} решения нерегулярных совместных систем (\ref{(36)})
с отличными от нуля $z_k$. Надо заменить $S_l$ на
$S_l+\sum_{k=1}^{s}{\varepsilon}_k^{l}$, где $s$
--- наименьшее натуральное число, при котором соответствующий
производящий многочлен отличен от тождественного нуля. Затем,
вычислив предел (\ref{(41)}), найти  нерегулярный производящий
многочлен $G_n^*(z)$ и его ненулевые корни $z_k$.

\subsection{Регуляризация. Второй способ} \label{s4.2}

Приведем некоторые условия $({\mathcal B}')$-разрешимости
переопределенной системы (\ref{(36)}) (а значит, для допустимого
решения задачи (\ref{(2)})) в случае $n=s\mu-1$.

\begin{theorem} \label{theorem5} При $n=s\mu-1$ для $({\mathcal
B}')$-разрешимости переопределенной системы $(\ref{(36)})$
необходимо $r_s(\omega)=0$.
 \end{theorem}

Доказательство следует непосредственно из леммы~\ref{lemma5} и теорем~\ref{theorem2} и~\ref{theorem4}.
По-видимому это необходимое условие является и достаточным. Мы
докажем это в случае $\mu=2$.

\begin{lemma} \label{lemma8} В условиях теоремы $\ref{theorem2}$
 имеем $ g_{\mu-1}(\omega)=(-1)^{s+1} r_{s}(\omega)$
 при $\mu=2$, т.е.
\begin{equation}
G_{n}(\omega,2;\,z)=(-1)^{s+1} r_{s}(\omega)\cdot z \cdot
{\mathcal R}_{s-1}(\omega;z^{2}),\qquad n=2s-1.
 \label{(42)}
\end{equation}
\end{lemma}

\begin{Proof} Минор $M_n$ (размера $n\times n$, $n=2s-1$),
дополнительный к $z^1$ в (\ref{(17)}), имеет структуру, позволяющую
упростить его путем определенного <<прореживания>>.

На шаге $k=1$ строка минора $M_n$ с номером $n-1$ содержит
единственный ненулевой элемент с номером $(n-1,3)$, равный
единице. Через $M_{n-1}$ обозначим дополнительный к ней минор
(размера $(n-1)\times (n-1)$). Далее выполняются сходные операции.
На шаге $k\ge 2$ строка $(n-2k+1)$ минора $M_{n-k+1}$ содержит
единственный ненулевой элемент с номером $(n-2k+1,2+k)$, равный
единице. Через $M_{n-k}$ обозначим дополнительный к ней минор. На
последнем шаге $k=s-1$ строка минора $M_{s+1}$ с номером $2$
содержит единственный ненулевой элемент с номером $(2,s+1)$,
равный единице. Через $M_{s}$ обозначим дополнительный к ней
минор. Учтем знаки, с которыми берутся миноры в указанном
разложении и, просуммировав номера
$1+\sum_{k=1}^{s-1}((n-2k+1)+(2+k))$, получим
$g_1(\omega)=-(-1)^{s(s-1)/2} M_s$.

Заметим, что полученный определитель $M_s$ зеркально симметричен к
трехдиагональному определителю $J_s$ Якоби с постоянными
элементами $\omega$ на главной диагонали и единицами на двух
прилегающих диагоналях (т.е. элементы $\alpha_{k,l}$ определителя
$J_s$ имеют вид $\alpha_{k,k}=\omega$, $\alpha_{k,k+1}=1$,
$\alpha_{k-1,k}=1$, остальные элементы нулевые). Очевидно,
$M_s=(-1)^{s(s-1)/2} J_s$. Для $J_k$ хорошо известна рекуррентная
формула (см., напр.,\cite{PRASOLOV})
$$
J_k=\omega J_{k-1}-J_{k-2},\qquad J_{-1}=0,\quad J_0=1,\quad
k=\overline{1,s}.
$$
Учитывая (\ref{(22)}), (\ref{(25)}), видим, что
$J_s(\omega)=r_s(-\omega)=(-1)^sr_s(\omega)$. \end{Proof}

 \begin{theorem} \label{theorem6} При $\mu=2$, $n=s\mu-1$ $(s\ge 2)$  для
$({\mathcal B}')$-разрешимости системы $(\ref{(36)})$ необходимо и
достаточно условие $r_{s}(\omega)=0$ $($ т.е. $\omega\in
\Omega_s)$. При этом $m=n-1$, а решение $\{z_k\}$ состоит из
корней многочлена ${\mathcal R}_{s-1}(\omega;z^{2})$.
 \end{theorem}

\begin{Proof} Необходимость условия $r_{s}(\omega)=0$
вытекает из леммы~\ref{lemma8} и теоремы~\ref{theorem5}. Докажем достаточность. Пусть
$\omega$  --- корень многочлена $r_s$, тогда $r_{s-1}(\omega)\ne
0$ (учитываем ортогональность многочленов) и из (\ref{(22)}) получаем
\begin{equation}
\omega=-r_{s-2}(\omega)/r_{s-1}(\omega),\qquad r_{s-1}(\omega)\ne
0.
 \label{(43)}
\end{equation}
Рассмотрим варьированную задачу (\ref{(16)}), заменив $\omega$ на
$\omega+\varepsilon$, $\varepsilon>0$. В силу (\ref{(43)}) и
леммы~\ref{lemma6} при $\varepsilon\to 0$ корни многочлена
${\mathcal R}_{s-1}(\omega+\varepsilon;z^{2})$ стремятся к корням
$z_j$ ($j=\overline{2,2s-1}$) многочлена ${\mathcal
R}_{s-1}(\omega;z^{2})$. Корни $z_j$ попарно различны и лежат на
единичной окружности (они составляют множество (\ref{(30)}), где
следует заменить $t_k$ на $t_k^2$). Поэтому многочлен
$G_{n}(\omega+\varepsilon,2;\,z)$ при достаточно малых
$\varepsilon$ имеет $n=2s-1$ различных корней $\xi_k(\varepsilon)$
$(k=\overline{1,n})$, один из которых, пусть $\xi_1$, равен нулю
и, как уже сказано, $\xi_{k}(\varepsilon)\to z_k$
$(k=\overline{2,n})$ при $\varepsilon\to 0$. В частности, из
(\ref{(30)}) следует, что
$|\xi_k(\varepsilon)-\xi_m(\varepsilon)|\ge {\rm const}>0$, $k\ne
m$, при достаточно малых $\varepsilon$.

Следовательно, получается $({\mathcal B}'')$-регулярная система
(\ref{(16)}). Совершая предельный переход при $\varepsilon\to 0$,
видим, что $Z_k(\varepsilon)\to Z_k=Z_k(0)$ с некоторыми конечными
$Z_k(0)$, поскольку определитель первых $n$ уравнений системы
(\ref{(16)}) отделен от нуля и стремится к отличной от нуля
константе.

Остается показать, что $Z_1=0$. Воспользуемся формулой (\ref{(6)}). Тогда
${Z}_1=({\mathcal{L}}_1\cdot {\mathcal{S}})$, где $
{\mathcal{S}}=(S_0,\ldots,S_{n-1})$ и
$S_{n-1}=\omega+\varepsilon$, $S_{n-3}=1$ (остальные равны нулю).
Положим
$$
\rho_2^{(1)}(\{\xi_k\})=\sum_{2 \le j_1<j_2 \le
n}{{\xi}_{j_1}{\xi}_{j_2}},\qquad\rho_2^{(1)}(\{z_k\})=\sum_{2 \le
j_1<j_2 \le n}{{z}_{j_1}{z}_{j_2}}.
$$
Из (\ref{(43)}) имеем $\rho_2^{(1)}(\{z_k\})=-\omega$, откуда с учетом (\ref{(6)})
находим
$$
Z_1(\varepsilon)=\frac{\rho_2^{(1)}(\{\xi_k\})+(\omega+\varepsilon)}
{\prod_{k=2}^{n}(\xi_1-\xi_k)}=
\frac{\varepsilon+(\rho_2^{(1)}(\{\xi_k\})-\rho_2^{(1)}(\{z_k\}))}
{\prod_{k=2}^{n}(\xi_1-\xi_k) }\to 0.
$$
Итак, варьированная система (\ref{(16)}) в пределе дает
$({\mathcal B}')$-разрешимую систему (\ref{(36)}) с $m=n-1$.
\end{Proof}

\subsection{Регуляризация. Третий способ} \label{s4.3}

Как и в предыдущем пункте, рассматриваем случай $n=s\mu -1$. Пусть
система (\ref{(16)}) совместна с некоторым $m<n$ и удовлетворяет
условию $({\mathcal B}')$. Из теоремы~\ref{theorem2} следует, что
в случае $n=s\mu -1$ при $s\ge 2$ и $\mu\ge 2$ система
(\ref{(16)}) всегда нерегулярна. Проведем ее регуляризацию,
заменив $S_{2n-1}=\sigma_n=0$ на $\sigma_n=\varepsilon$ с
некоторым параметром $\varepsilon\in{\mathbb C}$ (остальные
параметры не изменяются). В определителе (\ref{(17)}) при этом
изменится только один элемент в правом нижнем углу. Учитывая
(\ref{(35)}) и то, что новый производящий многочлен
$\hat{G}(\varepsilon,\omega,\mu;z)$ зависит от $\varepsilon$
линейно, получим
\begin{equation}
\hat{G}(\varepsilon,\omega,\mu;z)=\varepsilon {\mathcal
T}(\omega,\mu;z)+g_{{\mu}-1}(\omega) z^{{\mu}-1} \cdot {\mathcal
R}_{s-1}(\omega;z^{\mu}).
  \label{(44)}
\end{equation}

Найдем явный вид многочлена ${\mathcal T}$.
\medskip

\begin{lemma} \label{lemma9} Верна формула:
\begin{equation}
{\mathcal T}(\omega,\mu;z)=A(\omega) z^{(s-1)\mu} {\mathcal
R}_{s-1}(\omega;z^{-\mu}),
 \label{(45)}
\end{equation}
где $A(\omega)$ --- некоторый многочлен степени $n-1$, отличный от
тождественного нуля. В частности, при $A(\omega)\ne 0$ многочлен
${\mathcal T}$ имеет корни, обратные по отношению к корням
многочлена ${\mathcal R}_{s-1}(\omega;z^{\mu})$.
\end{lemma}

\begin{Proof} Многочлен ${\mathcal T}$ является
алгебраическим дополнением в определителе (\ref{(17)}) к элементу,
находящемуся в правом нижнем углу. Его матрицу обозначим через
$B=(b_{j,k})$ $(k=\overline{1,n},\;j=\overline{1,n})$. Построим
далее матрицу $A=(a_{j,k})$, у которой столбцы и строки зеркально
симметричны (т.е. расположены в обратном порядке) по отношению к
столбцам и строкам матрицы $B$ (поворот матрицы $B$ на $180^{0}$;
точнее, $a_{j,k}=b_{n+1-j,n+1-k}$). Заметим, что если перенести
последнюю строку (со степенями $z^j$) матрицы $A$ на первое место
и формально заменить в ней $z^k$ на $z^{n-1-k}$
$(k=\overline{0,n-1})$, то получится матрица той же структуры, что
в (\ref{(17)}), но на единицу меньшего размера (при условии
$({\mathcal B}')$ на $\sigma_k$). Это означет, что $\det
A=(-1)^{n-1} z^{n-1} G_{n-1}(\omega,\mu;1/z)$.

Как и в лемме~\ref{lemma7} доказывается, что коэффициенты $\tilde g_{m+k\mu}$
многочлена $G_{n-1}(\omega,\mu;\xi)$, $\xi=1/z$, равны нулю при
$m=\overline{0,\mu-3}$, $k\in {\mathbb N}$. (В этом случае линейно
зависимы столбцы с номерами $k\mu -1$, $k=\overline{1,s}$.) Отсюда
и из равенств (\ref{(23)}) аналогично (\ref{(35)}) получаем $\tilde
g_{{\mu}k-2}=\tilde g_{\mu-2} r_{k-1}$ при $k=\overline{1,s}$ и
$$
\det A=(-1)^{n-1} {\xi}^{1-n} G_{n-1}(\omega,\mu;\xi)=(-1)^{n-1}
  \tilde g_{{\mu}-2} {\xi}^{1-n} {\xi}^{{\mu}-2} \cdot {\mathcal
R}_{s-1}(\omega;{\xi}^{\mu}),
$$
где  $\tilde g_{{\mu}-2}(\omega)$ --- многочлен, не равный
тождественно нулю (поскольку при $\omega\to\infty$ старшее
слагаемое в $G_{n-1}(\omega,\mu;\xi)$ имеет вид
$\pm\omega^{n-1}(1+o(1)){\xi}^{n-1}$). Учитывая, что ${\mathcal
T}$ отличается от $\det A$ только знаком, получаем нужное.
\end{Proof}

В следующей теореме можно считать, что $\mu\ge 3$, поскольку при
$\mu=1,2$ полное решение задачи (\ref{(2)}) дано в
теоремах~\ref{theorem3} и \ref{theorem6}.

\begin{theorem} \label{theorem7} Пусть $n=s\mu -1$ $(s\ge 2,\;\mu\ge 2)$ и
$\omega\in \Omega_s$. Тогда система $(\ref{(36)})$ с исключенным
последним уравнением имеет допустимое решение, причем $m=
n-(\mu-1)$. Различные корни $z_k$, $k=\overline{1,m}$, являются
нулями многочлена ${\mathcal R}_{s-1}(\omega;z^{\mu})$ и при
$\omega=-2\cos \varphi_{\alpha}=2\cos \varphi_{s-\alpha+1}$
составляют множество
\begin{equation}
\{z_k\}=\{\sqrt[(s+1)\mu]{(-1)^{\alpha}}\}\setminus
\{\sqrt[\mu]{e^{i\varphi_{\alpha}}},
\sqrt[\mu]{e^{-i\varphi_{\alpha}}}\},\qquad \varphi_{\alpha}
 =\frac{\pi {\alpha}}{s+1},\quad {\alpha}=\overline{1,s}.
   \label{(46)}
\end{equation}
При этом
\begin{equation}
Z_k=-\frac{\sum_{j\ne k}z_j}{\prod_{j\ne k}(z_k-z_j)}=\frac{z_k
}{\prod_{j\ne k}(z_k-z_j)},\qquad k=\overline{1,m}.
      \label{(47)}
\end{equation}
В этом случае допустимое решение $($с условием $(\mathcal A)$$)$
имеет вид
$$
X_k=Z_k z_k^{n-1},\qquad e^{-i\lambda_k}=z_k,\qquad
k=\overline{1,m},
$$
а соответствующий АФО имеет порядок $m=n-(\mu-1)$ и решает задачу
$(\ref{(2)})$ на многочленах порядка $n-1$.
\end{theorem}

\begin{remark} \label{remark2} Аналогично примечанию~\ref{remark1a} отметим, что в силу
равенства $z_k^{{\tau}}=(-1)^{\alpha}$,
${\tau}:=(s+1)\mu=n+\mu+1$, нулевые степенные суммы $\sum_k
X_kz_k^{\beta}$ в (\ref{(4)}) порядка $\beta$ остаются нулевыми при замене
$\beta$ на $\beta+{\tau}\,k$, $k\in{\mathbb N}$. Следовательно,
отличными от нуля могут быть лишь степенные суммы порядка
\begin{equation}
\beta=\mu+{\tau}\,k\qquad \hbox{и}\qquad n+{\tau}\,k\le \beta\le
\tau+{\tau}k,\quad {\tau}=n+\mu+1.
 \label{(48)}
\end{equation}
Поэтому формула (\ref{(2)}) при $\mu\ge 2$ верна и для сходящихся
тригонометрических рядов, в которых отсутствуют гармоники с
номерами (\ref{(48)}).
\end{remark}

\begin{Proof} Пусть $\omega\in \Omega_s$. Тогда
$r_{s}(\omega)=0$, $r_{s-1}(\omega)\ne 0$ (учитываем
ортогональность многочленов). Выберем последовательность
$\{\omega_k\}$, $\omega_k\to \omega$ так, чтобы
$g_{{\mu}-1}(\omega_k)\ne 0$, $A(\omega_k)\ne 0$,
$r_{s-1}(\omega_k)\ne 0$. По лемме~\ref{lemma6} корни многочлена
${\mathcal R}_{s-1}(\omega_k;{z}^{\mu})$ сходятся к корням
многочлена ${\mathcal R}_{s-1}(\omega;{z}^{\mu})$ (попарно
различным и лежащим на единичной окружности симметрично
относительно действительной оси, см. (\ref{(30)})). Из
(\ref{(45)}) следует, что  корни многочлена ${\mathcal
T}(\omega,\mu;z)$ также сходятся к корням многочлена ${\mathcal
R}_{s-1}(\omega;{z}^{\mu})$. Далее, положим
$\alpha_k=g_{{\mu}-1}(\omega_k) r_{s-1}(\omega_k)$,
$N_k(z)={\mathcal R}_{s-1}(\omega_k;z^{\mu})/r_{s-1}(\omega_k)$,
так что (см. (\ref{(44)}))
$$
\hat{G}_k(z):=\hat{G}(\varepsilon_k,\omega_k,\mu;z)= \varepsilon_k
{\mathcal T}(\omega_k,\mu;z)+\alpha_k z^{{\mu}-1} \cdot N_k(z),
$$
причем параметр $\varepsilon_k$ выберем из условия:
$$
\varepsilon_k {\mathcal T}(\omega_k,\mu;z)=-2^{-1}\alpha_kM_k(z),
$$
где $M_k$--- унитарные многочлены (с единичными старшими
коэффициентами) степени $\mu(s-1)$. Тогда
\begin{equation}
\hat{G}_k(z)=\alpha_k (z^{{\mu}-1} -2^{-1})N_k(z)+\alpha_k
2^{-1}(N_k(z)-M_k(z)),
    \label{(49)}
\end{equation}
где, как уже сказано, $M_k$, $N_k$--- унитарные многочлены степени
$\mu(s-1)$, корни которых при $k\to\infty$ сходятся к корням
многочлена ${\mathcal R}_{s-1}(\omega;{z}^{\mu})$.

Поскольку $\alpha_k$ ограничены, то, очевидно,
$\alpha_k(N_k(z)-M_k(z))$ равномерно стремится к нулю на любом
компактном множестве. Поэтому по теореме Руше многочлены
$\hat{G}_k$ при достаточно больших $k\ge k_0$ имеют $n=s\mu-1$
различных корней, которые сходятся к $\sqrt[\mu-1]{1/2}$ и корням
многочлена ${\mathcal R}_{s-1}(\omega;{z}^{\mu})$. Пусть для
определенности корни $z^{(k)}_1,\ldots,z^{(k)}_{\mu-1}$ сходятся к
$\sqrt[\mu-1]{1/2}$, а остальные корни
$z^{(k)}_{\mu},\ldots,z^{(k)}_{n}$ сходятся к корням
$z_{\mu},\ldots,z_{n}$ многочлена ${\mathcal
R}_{s-1}(\omega;{z}^{\mu})$.

Получается последовательнось $({\mathcal B}'')$-регулярных систем
(\ref{(16)}) с производящими многочленами $\hat{G}_k$ и решениями
$\{Z^{(k)}_j\}$, $\{z^{(k)}_j\}$. Покажем, что $Z^{(k)}_j\to 0$
при $k\to\infty$ для номеров $j=\overline{1,\mu-1}$.
Действительно, возьмем, к примеру, $j=1$. Первые $n$ уравнений
варьрованной системы (\ref{(16)}) имеют вид (\ref{(4)}) с
условиями $({\mathcal B}')$, где в соответствии с нумерацией
уравнений (\ref{(16)}) следует заменить $\mu$ на $n-\mu$. По
формуле (\ref{(6)}) получаем
 $$
 Z_1^{(k)}=
 \frac{\omega_k+(-1)^{\mu}\rho_{\mu}(z^{(k)}_2,\ldots,z^{(k)}_n)}
 {\prod_{j=2}^{n}
(z^{(k)}_1-z^{(k)}_j)} = \frac{b_k+ a_k} {\prod_{j=2}^{n}
(z^{(k)}_1-z^{(k)}_j)},
$$
где
$$
b_k=\omega+ (-1)^{\mu}\rho_{\mu}(z^{(k)}_{\mu},\ldots,z^{(k)}_n),
$$
$$
a_k:=(\omega_k-\omega)+
(-1)^{\mu}\sum_{m=1}^{\mu-2}\rho_{m}(z^{(k)}_{2},\ldots,z^{(k)}_{\mu-1})
 \rho_{\mu-m}(z^{(k)}_{\mu},\ldots,z^{(k)}_{n}).
$$
При $\mu=2$ последняя сумма считается равной нулю. По теореме
Виета для многочлена ${\mathcal R}_{s-1}(\omega;{z}^{\mu})$ имеем
(см. (\ref{(43)}))
$$
\omega+ (-1)^{\mu}\rho_{\mu}(z_{\mu},\ldots,z_n)=0,\qquad
\rho_{\mu-m}(z_{\mu},\ldots,z_{n})=0,\qquad m=\overline{1,\mu-2}.
$$
Поэтому, учитывая, что $z^{(k)}_{j}\to z_j$ и
$|z^{(k)}_{1}-z^{(k)}_{j}|\ge {\rm const}>0$, $j\ne 1$, получаем
$b_k\to 0$, $a_k\to 0$, $Z_1^{(k)}\to 0$ ($k\to \infty$).

По тем же формулам (\ref{(6)}) получается и сходимость $Z_j^{(k)}$
к некоторым конечным $Z_{j+1-\mu}$, $j=\overline{\mu,n}$ (индексы
изменены для удобства). Для определения $Z_j$ остается решить
первые $m=n-(\mu-1)$ уравнений предельной системы вида
(\ref{(36)}). Она имеет вид (\ref{(4)}), где столбец сводных
членов содержит единственный ненулевой элемент (равный единице),
находящийся на предпоследнем месте с номером $m-1$. По формуле
(\ref{(6)}) получается первое равенство в (\ref{(47)}), а во
втором учитывается равенство $\sum_{j=1}^{m}z_j=0$ (многочлен
${\mathcal R}_{s-1}(\omega;{z}^{\mu})$ не содержит степени
$z^{n-\mu}$).

Итак, в пределе получаем тождества (\ref{(36)}) при $l=\overline{0,2n-2}$
и $m=n+1-\mu$. Последнее равенство в (\ref{(36)}) приходится исключить,
поскольку из приведенного доказательства предел
$\lim_{k\to\infty}\varepsilon_k$ не определяется (он существует в
силу существования пределов от $Z_j^{(k)}$ и $z_j^{(k)}$ и, скорее
всего, равен нулю). В результате приходим к (\ref{(15)}), а значит и к (\ref{(4)})
без последних уравнений. Вещественность АЧС следует из леммы~\ref{lemma1}
(исключение последнего уравнения в лемме роли не играет, т.к.
$m<n$). \end{Proof}

\bigskip

\end{document}